\documentclass[10pt,twoside]{article}
\usepackage{graphicx}
\usepackage{amsmath}
\usepackage{Latex-document}

\newcommand{\co}{{\colon\thinspace}}
\newcommand{\Zed}{{\bf{Z}}}
\newcommand{\Real}{{\bf{R}}}
\newcommand{\Complex}{{\bf{C}}}

\title{\bf Representations of Braid Groups\vskip 6mm}

\author{S. Bigelow\vspace*{-0.5cm}\thanks{Department of Mathematics \& Statistics, University of Melbourne,
Victoria 3010, Australia. E-mail: bigelow@unimelb.edu.au}}
\date{\vspace{-8mm}}

\begin{document}
\maketitle

\thispagestyle{first} \setcounter{page}{37}

\begin{abstract}\vskip 3mm
In this paper we survey some work on
representations of $B_n$
given by the induced action on a homology module of some space.
One of these, called the Lawrence-Krammer representation,
recently came to prominence when it was shown to be faithful for all $n$.
We will outline the methods used,
applying them to a closely related representation
for which the proof is slightly easier.
The main tool is the Blanchfield pairing,
a sesquilinear pairing between elements of relative homology.
We discuss two other applications of the Blanchfield pairing,
namely a proof that the Burau representation is not faithful for large $n$,
and a homological definition of the Jones polynomial.
Finally,
we discuss possible applications
to the representation theory of the Hecke algebra,
and ultimately of the symmetric group over fields of non-zero characteristic.

\vskip 4.5mm


\end{abstract}

\vskip 12mm

\section{Introduction} \setzero
\vskip-5mm \hspace{5mm }

Artin's braid group $B_n$ was originally defined as a group
of geometric braids in $\Real^3$.
Representations of $B_n$ have been studied
for their own intrinsic interest,
and also in connection to other areas of mathematics,
most notably to knot invariants such as the Jones polynomial.

We will use the definition of $B_n$
as the mapping class group of
an $n$-times punctured disk $D_n$.
A rich source of representations of $B_n$
is the induced action
on homology modules of spaces related to $D_n$.
The Burau representation,
one of the simplest and best known representations of braid groups,
is most naturally defined as the induced action of $B_n$
on the first homology of a cyclic covering space of $D_n$.
Lawrence \cite{rL90}
extended this idea to configuration spaces in $D_n$,
and was able to obtain
all of the so-called Temperley-Lieb representations.

Lawrence's work seems to have received very little attention
until one of her homological representations
was shown to be faithful,
thus proving that braid groups are linear
and solving a longstanding open problem.
Two independent and very different proofs of this have appeared in
\cite{sB01} and \cite{dK02}.
In this paper we will outline the former,
emphasising the importance of the Blanchfield pairing.
We then discuss two other applications of the Blanchfield pairing,
namely
the proof that the Burau representation is not faithful for large $n$,
and a homological definition of the Jones polynomial of a knot.
We conclude with some speculation
on possible future applications to representations of Hecke algebras
when $q$ is a root of unity.
These are related to representations of the symmetric group $S_n$
over fields of bad characteristic,
that is, fields in which $n!=0$.

%--------------------------------------------------
\section{The Lawrence-Krammer representation}
\vskip-5mm \hspace{5mm }

Let $D$ be the unit disk centred at the origin in the complex plane.
Fix arbitrary real numbers $-1 < p_1 < \dots < p_n < 1$,
which we will call ``puncture points''.
Let
$$D_n = D \setminus \{p_1,\dots,p_n\}$$
be the $n$-times punctured disk.
The braid group $B_n$ is the mapping class group of $D_n$,
that is,
the set of homeomorphisms from $D_n$ to itself
that act as the identity on $\partial D$,
taken up to isotopy relative to $\partial D$.
Let $C$ be
the space of all unordered pairs of distinct points in $D_n$.

Suppose $x$ is a point in $D_n$,
and $\alpha$ is a simple closed curve in $D_n$
enclosing one puncture point and not enclosing $x$.
Let $\gamma \co I \to C$ be the loop in $C$ given by
$$\gamma(s) = \{x,\alpha(s)\}.$$
Now suppose $\tau_1$ and $\tau_2$ are paths in $D_n$
such that $\tau_1 \tau_2$ is a simple closed curve
that does not enclose any puncture points.
Let $\tau \co I \to C$ be the loop in $C$ given by
$$\tau(s) = \{\tau_1(s),\tau_2(s)\}.$$
Let
$$\Phi \co \pi_1(C) \to \langle q \rangle \oplus \langle t \rangle$$
be the unique homomorphism such that
$\Phi(\gamma) = q$ and $\Phi(\tau) = t$
for any $\gamma$ and $\tau$ defined as above.
For a proof of the existence and uniqueness of such a homomorphism,
see \cite{PP02}.
Let $\tilde{C}$ be the connected covering space of $C$
whose fundamental group is the kernel of $\Phi$.

The second homology $H_2(\tilde{C})$
is a module over $\Zed[q^{\pm 1},t^{\pm 1}]$,
where $q$ and $t$ act by covering transformations.
The Lawrence-Krammer representation of $B_n$
is the induced action of $B_n$ on $H_2(\tilde{C})$
by $\Zed[q^{\pm 1},t^{\pm 1}]$-module automorphisms.
More precisely,
given an element of $B_n$
represented by a homeomorphism $\sigma \co D_n \to D_n$,
consider the induced action $\sigma \co C \to C$.
There is a unique lift
$\tilde{\sigma} \co \tilde{C} \to \tilde{C}$
that acts as the identity on $\partial\tilde{C}$.
This induces an automorphism of $H_2(\tilde{C})$,
which can be shown to respect the $\Zed[q^{\pm 1},t^{\pm 1}]$-module structure.
See \cite{PP02} for details.

%--------------------------------------------------
\section{The Blanchfield pairing}
\vskip-5mm \hspace{5mm }

Let $\epsilon > 0$ be small.
We define $P,B \subset C$ as follows.
Suppose $\{x,y\}$ is a point in $C$.
We say $\{x,y\} \in P$ if either $|x-y| \le \epsilon$,
or there is a puncture point $p_i$ such that
$|x-p_i| \le \epsilon$ or $|y-p_i| \le \epsilon$.
We say $\{x,y\} \in B$ if $x \in \partial D$ or $y \in \partial D$.

For $u \in H_2(\tilde{C},\tilde{P})$ and $v \in H_2(\tilde{C},\tilde{B})$
let $(u \cdot v) \in \Zed$ denote the standard algebraic intersection number.
We define an intersection pairing
$$\langle \cdot,\cdot \rangle
  \co H_2(\tilde{C},\tilde{P}) \times H_2(\tilde{C},\tilde{B})
  \to \Zed[q^{\pm 1},t^{\pm 1}]$$
by
$$\langle u,v \rangle
  = \sum_{i,j \in \Zed} (u \cdot q^i t^j v)q^i t^j.$$
For a proof that these are well-defined, see \cite[Appendix E]{aK96},
where the following properties are also proved.

For $u \in H_2(\tilde{C},\tilde{P})$,
$v \in H_2(\tilde{C},\tilde{B})$,
$\sigma \in B_n$, and $\lambda \in \Zed[q^{\pm 1},t^{\pm 1}]$, we have
$$\langle \sigma u,\sigma v \rangle = \langle u,v \rangle,$$
and
$$\langle \lambda u,v \rangle = \lambda \langle u,v \rangle
                              = \langle u,\bar{\lambda}v \rangle,$$
where $\bar{\lambda}$ is the image of $\lambda$
under the automorphism of $\Zed[q^{\pm 1},t^{\pm 1}]$
taking $q$ to $q^{-1}$ and $t$ to $t^{-1}$.

%--------------------------------------------------
\section{A faithful representation}
\vskip-5mm \hspace{5mm }

The aim of this section is to outline a proof of the following.

{\bf Theorem.}
\it
Let $\tilde{C}$ and $\tilde{P}$ be as above.
The induced action of $B_n$ on $H_2(\tilde{C},\tilde{P})$ is faithful.
\rm

For the details,
the reader is referred to \cite{sB01},
where the same techniques are used
to prove that $B_n$ acts faithfully on $H_2(\tilde{C})$.
Our use of relative homology here
actually simplifies the argument somewhat.

There is a slight technical difficulty in defining
the action of $B_n$ on $H_2(\tilde{C},\tilde{P})$.
Namely, the action of a braid on $C$
need not preserve the set $P$.
Thus we should really take a limit as $\epsilon$ approaches $0$.
The representation obtained is very similar to
the Lawrence-Krammer representation,
but has a slightly different module structure,
as discussed in \cite{sBgeorgia}.

Let $E$ be the straight edge from $p_1$ to $p_2$.
Let $E'$ be the set of points in $C$
of the form $\{x,y\}$, where $x,y \in E$.
Let $\tilde{E}'$ be a lift of $E'$ to $\tilde{C}$.
This represents an element of $H_2(\tilde{C},\tilde{P})$,
which we will call $e$.
Let $F_1$ and $F_2$ be parallel vertical edges
with endpoints on $\partial D$,
passing between $p_2$ and $p_3$.
Let $F'$ be the set of points in $C$
of the form $\{x,y\}$, where $x \in F_1$ and $y \in F_2$.
Let $\tilde{F'}$ be a lift of $F'$ to $\tilde{C}$.
This represents an element of $H_2(\tilde{C},\tilde{B})$,
which we will call $f$.
Note that
$$\langle e,f \rangle = 0,$$
since $E'$ and $F'$ are disjoint in $C$.

Suppose the action of $B_n$ on $H_2(\tilde{C},\tilde{P})$ is not faithful.
It is not hard to show that there must be a braid $\sigma$
in the kernel of this representation
such that $\sigma(E)$ is not isotopic to $E$ relative to endpoints.
Now $\sigma(e) = e$, so
$$\langle \sigma(e),f \rangle = 0.$$
From this, we will derive a contradiction.

By applying an isotopy,
we can assume $\sigma(E)$ intersects $F_1$ and $F_2$ transversely
at a minimal number of points $k > 0$.
Let $x_1,\dots,x_k$ be the points of $\sigma(E) \cap F_1$,
and let $y_1,\dots,y_k$ be the points of $\sigma(E) \cap F_2$,
numbered from top to bottom in both cases.

For $i,j \in \{1,\dots,k\}$,
let $a_{i,j}$ and $b_{i,j}$ be the unique integers such that
$\sigma(\tilde{E}')$ intersects $q^{a_{i,j}}t^{b_{i,j}} \tilde{F'}$
at a point in the fibre over $\{x_i,y_j\}$,
and let $\epsilon_{i,j}$ be the sign of that intersection.
Then
$$\langle \sigma(e),f \rangle
  = \sum_{i=1}^k \sum_{j=1}^k \epsilon_{i,j} q^{a_{i,j}} t^{b_{i,j}}.$$

To calculate $a_{i,j}$ and $b_{i,j}$,
it is necessary to specify choices of lift for $E'$ and $F'$.
We will not do this
since we only need to calculate differences
$a_{i',j'} - a_{i,j}$ and $b_{i',j'} - b_{i,j}$.
To do this,
let $\gamma$ be a path in $C$
that goes from $\{x_i,y_j\}$ to $\{x_{i'},y_{j'}\}$ in $\sigma(E')$,
and then back to $\{x_i,y_j\}$ in $F'$.
Then
$$q^{(a_{i',j'} - a_{i,j})} t^{(b_{i',j'} - b_{i,j})} = \Phi(\gamma).$$
From this we can prove the following.

{\bf Lemma.}
\it
For all $i,j \in \{1,\dots,k\}$ we have
\begin{itemize}
\item $a_{i,j} = \frac{1}{2} (a_{i,i} + a_{j,j})$,
\item if $b_{i,j} > b_{i,i}$
      then $a_{i,j'} > a_{i,i}$ for some $j' = 1,\dots,k$,
\item if $b_{i,j} > b_{j,j}$
      then $a_{i',j} > a_{j,j}$ for some $i' = 1,\dots,k$.
\end{itemize}
\rm

The first of these is \cite[Lemma 2.1]{sB01},
and the second and third follow from the proof of \cite[Claim 3.4]{sB01}.
We now sketch the proof of the second
in the special case where
$y_i$ lies between $x_i$ and $y_j$ along $\sigma(E)$.

Let $\alpha$ be the path from $y_i$ to $y_j$ along $\sigma(E)$.
Let $\beta$ be the path from $y_j$ to $y_i$ along $F_2$.
Then $b_{i,j} - b_{i,i}$ is two times
the winding number of $\alpha\beta$ around $x_i$.
In particular, this winding number is positive.

Let $D_1$ be the once punctured disk $D \setminus \{x_i\}$,
and let $\tilde{D}_1$ be its universal cover.
Let $\tilde{\alpha}\tilde{\beta}$ be a lift of $\alpha\beta$ to $D_1$.
This is a path from a point in the fibre over $y_i$
to a ``higher'' point in the fibre over $y_i$.

Let $F_2^+$ be the the segment of $F_2$
going from $y_i$ upwards to $\partial D$.
Let $\tilde{F}_2^+$ be
the lift of $F_2^+$ to $\tilde{D}_1$
that has an endpoint at $\tilde{\alpha}(0)$.
In order to reach a higher sheet in $\tilde{D}_1$,
$\tilde{\alpha}$ must intersect $\tilde{F}_2^+$.
Let $\tilde{\gamma}$ be the loop in $\tilde{D}_1$
that follows $\tilde{\alpha}$
to the first point of intersection with $\tilde{F}_2^+$,
and then follows $\tilde{F}_2^+$ back to $\tilde{\alpha}(0)$.

Let $\gamma$ be the projection of $\tilde{\gamma}$ to $D_1$.
This travels along $\sigma(E)$
from $y_i$ to some point $y_{j'} \in F_2$,
then along $F_2$ back to $y_i$.
Then $a_{i,j'} - a_{i,i}$ is the total winding number of $\gamma$
around the puncture points.
We must show that this winding number is positive.

By construction,
$\tilde{\gamma}$ is a simple closed curve in $\tilde{D}_1$.
By the Jordan curve theorem, it must bound a disk $\tilde{B}$.
Let $B$ be the projection of $\tilde{B}$ to $D_1$.
This is an immersed disk in $D$, whose boundary is $\gamma$.
Note that $\tilde{\gamma}$ passes anticlockwise around $\tilde{B}$,
since the puncture $x_i$ lies to its right.
Thus the total winding number of $\gamma$ around the puncture points
is equal to the total number of puncture points contained in $B$,
counted with multiplicities.

It remains to show that $B$ intersects at least one puncture point.
Suppose not.
Then $B$ is an immersed disk in $D_n$.
Using an ``innermost disk'' argument,
one can find an embedded disk $B'$ in $D_n$
whose boundary consists of a subarc of $\sigma(E)$ and a subarc of $F_2$.
Using $B'$, one can isotope $\sigma(E)$
so as to have fewer points of intersection with $F_2$,
thus contradicting our assumptions.

This completes the proof of the second part of the lemma
in the case where
$y_i$ lies between $x_i$ and $y_j$ along $\sigma(E)$.
The remaining case, where $x_i$ lies between $y_i$ and $y_j$,
is only slightly trickier.
The third part of the lemma is similar to the second.
The first part of the lemma is much easier.

We now return to the proof of the theorem.
Let $a$ be the maximum of all $a_{i,j}$.
Let $b$ be the maximum of $\{b_{i,j} : a_{i,j} = a\}$.
Suppose $i,j \in \{1,\dots,k\}$ are such that
$a_{i,j} = a$ and $b_{i,j} = b$,
and also $i',j' \in \{1,\dots,k\}$ are such that
$a_{i',j'} = a$ and $b_{i',j'} = b$.
I claim that $\epsilon_{i,j} = \epsilon_{i',j'}$.
From this claim,
it follows that all occurrences of $q^a t^b$ in the expression
$$\langle \sigma(e),f \rangle
  = \sum_{i=1}^k \sum_{j=1}^k \epsilon_{i,j} q^{a_{i,j}} t^{b_{i,j}} $$
occur with the same sign,
so the coefficient of $q^a t^b$ is non-zero in $\langle \sigma(e),f \rangle$.
This provides our desired contradiction,
and completes the proof of the theorem.
It remains to prove that $\epsilon_{i',j'} = \epsilon_{i,j}$.

Using the above lemma,
it is not hard to show that
$a_{i,i} = a_{j,j} = a$ and $b_{i,i} = b_{j,j} = b$.
Similarly,
$a_{i',i'} = a_{j',j'} = a$ and $b_{i',i'} = b_{j',j'} = b$.
We will only need the equalities
$$b_{i,i} = b_{i,j} = b_{j,j} = b_{i',i'} = b_{i',j'} = b_{j',j'}.$$
In fact, we only need these modulo two.

Orient $\sigma(E)$ so that it crosses $F_1$ from left to right at $x_i$.
Let $\gamma$ be the path in $C$
which goes from $\{x_i,y_i\}$ to $\{x_{i'},y_{i'}\}$ in $E'$
and then back to $\{x_i,y_i\}$ in $F'$.
Now
$b_{i',i'} - b_{i,i}$ is the exponent of $t$ in $\Phi(\gamma)$.
The fact that this is an even number
means that the pair of points in $D_n$
do not ``switch places'' when they go around this loop.
Thus $\sigma(E)$ crosses $F_1$ from left to right at $x_i'$.
By similar arguments,
\begin{itemize}
\item $\beta(E)$ intersects $F_1$ with the same sign at $x_i$ and $x_{i'}$,
\item $\beta(E)$ intersects $F_2$ with the same sign at $y_j$ and $y_{j'}$,
\item $x_i$ occurs before $y_j$ and
      $x_{i'}$ occurs before $y_{j'}$
      with respect to the orientation of $\sigma(E)$.
\end{itemize}
It is now intuitively clear that
$E'$ must intersect $F'$ with the same signs at
$\{x_i,y_j\}$ and $\{x_{i'},y_{j'}\}$.
This can be proved rigorously
by careful consideration of orientations of these surfaces,
as discussed in \cite[Section 2.1]{sB01}.
It follows that $\epsilon_{i,j} = \epsilon_{i',j'}$,
which completes the proof of the theorem.

%--------------------------------------------------
\section{The Burau representation}
\vskip-5mm \hspace{5mm }

The proof that the Lawrence-Krammer representation is faithful
basically reduces to proving that
the Blanchfield pairing detects whether
corresponding edges in the disk can be isotoped to be disjoint.
A converse to this idea leads to a proof that
the Burau representation is not faithful for large $n$.

The Burau representation can be defined
by a similar but simpler construction
to that of the Lawrence-Krammer representation.
Let
$$\Phi \co \pi_1(D_n) \to \langle q \rangle$$
be the homomorphism that sends
each of the obvious generators to $q$.
Let $\tilde{D}_n$ be the corresponding covering space.
The Burau representation is the induced action of $B_n$
on $H_1(\tilde{D}_n)$ by $\Zed[q^{\pm 1}]$-module automorphisms.

Let $P$ be an $\epsilon$-neighbourhood of the puncture points,
and let $\tilde{P}$ be the preimage of $P$ in $\tilde{D}_n$.
The Blanchfield pairing in this context is a sesquilinear pairing
$$\langle \cdot,\cdot \rangle \co
  H_1(\tilde{D}_n,\tilde{P}) \times H_1(\tilde{D}_n,\partial\tilde{D}_n)
  \to \Zed[q^{\pm 1}].
$$

Let $E$ be the straight edge from $p_1$ to $p_2$.
Let $\tilde{E}$ be a lift of $E$ to $\tilde{D}_n$.
This represents an element of $H_1(\tilde{D}_n,\tilde{P})$,
which we will call $e$.
Let $F$ be a vertical edge
with endpoints on $\partial D$,
passing between $p_{n-1}$ and $p_n$.
Let $\tilde{F}$ be a lift of $F$ to $\tilde{D}_n$.
This represents an element of $H_1(\tilde{D}_n,\partial\tilde{D}_n)$,
which we will call $f$.
The following is \cite[Theorem 5.1]{sBkorea}.

{\bf Theorem.}
\it
Let $E$, $e$, $F$ and $f$ be as above.
The Burau representation of $B_n$ is unfaithful
if and only if there exists $\sigma \in B_n$
such that $\langle \sigma(e),f \rangle = 0$,
but $\sigma(E)$ is not isotopic relative to endpoints
to an edge that is disjoint from $F$.
\rm

Using this theorem,
one can show that the Burau representation of $B_n$ is not faithful
simply by exhibiting the required edges $\sigma(E)$ and $F$.
Such edges can be found by hand in the case $n=6$.
In the case $n=5$, they can be found by a computer search,
and then laboriously checked by hand.
The case $n=4$ seems to be beyond the reach of any known computer algorithm.
This is the last open case,
since the Burau representation is known to be faithful for $n \le 3$.

%--------------------------------------------------
\section{The Jones polynomial}
\vskip-5mm \hspace{5mm }

In this section,
we use the Blanchfield pairing to give
a homological definition of the Jones polynomial of a knot or link.
The Jones polynomial was defined in \cite{vJ85}
using certain algebraically defined representations of braid groups.
No satisfactory geometric definition is known,
but some insight might be offered
by defining the representations homologically
and using the Blanchfield pairing.
This was the original motivation for Lawrence
to study homological representations of braid groups.

A {\em geometric braid} $\sigma \in B_n$ is
a collection of $n$ disjoint edges in $\Complex \times \Real$
with endpoints $\{1,\dots,n\} \times \{0,1\}$,
such that each edge goes
from $\Complex \times \{0\}$ to $\Complex \times \{1\}$
with a constantly increasing $\Real$ component.
The correspondence between geometric braids
and elements of the mapping class group
is described in \cite{sB01},
and in many other introductory expositions on braids.
The {\em plat closure} of a geometric braid $\sigma \in B_{2n}$
is the knot or link obtained
by using straight edges to connect $(2j-1,k)$ to $(2j,k)$
for each $j=1,\dots,n$ and $k=0,1$.
Every knot or link can be obtained in this way.

Let $C$ be the configuration space of
unordered $n$-tuples of distinct points in $D_{2n}$.
Let
$$\Phi \co \pi_1(C) \to \langle q \rangle \oplus \langle t \rangle$$
be defined as in Section 2.
Namely,
if $\gamma$ is any loop in which one of the $n$-tuple
winds anticlockwise around a puncture point,
and $\tau$ is any loop in which two of the $n$-tuple
exchange places by an anticlockwise half twist,
then $\Phi(\gamma) = q$ and $\Phi(\tau) = t$.
Let $\tilde{C}$ be the covering space corresponding to $\Phi$.

Define $P,B \subset C$ similarly to those of Section 3.
The Blanchfield pairing is a sesquilinear pairing
$$\langle \cdot,\cdot \rangle
  \co H_2(\tilde{C},\tilde{P}) \times H_2(\tilde{C},\tilde{B})
  \to \Zed[q^{\pm 1},t^{\pm 1}].$$

For $k=1,\dots,n$,
let $F_k$ be the straight edge from $p_{2k-1}$ to $p_{2k}$.
Let $F$ be the set of points in $C$ of the form $\{x_1,\dots,x_n\}$
where $x_i \in F_i$.
Let $\tilde{F}$ be a lift of $F$ to $\tilde{C}$.
This represents an element of $H_n(\tilde{C},\tilde{P})$,
which we call $f$.
For $k=1,\dots,n$,
let $e_k \co S^1 \to D_n$ be a figure-eight
around $p_{2k-1}$ and $p_{2k}$
in a small neighbourhood of $F_k$.
Let $e$ be the map from the $n$-torus $(S^1)^n$ to $C$ given by
$$e(s_1,\dots,s_n) = \{e_1(s_1),\dots,e_n(s_n)\}.$$
Let $\tilde{e}$ be a lift of $e$ to $\tilde{C}$.
This represents an element of $H_n(\tilde{C})$,
and hence of $H_n(\tilde{C},\tilde{B})$,
which we also call $e$.

The main result of \cite{sBkyoto} is the following.

{\bf Theorem.}
\it
Let $e$ and $f$ be as above
and suppose $\sigma \in B_{2n}$.
The Jones polynomial of the plat closure of $\sigma$
is
$$\langle \sigma(e),f) \rangle |_{(t=-q^{-1})},$$
up to sign and multiplication by a power of $q^{\frac{1}{2}}$.
\rm

Here, the Jones polynomial is normalised so that
the Jones polynomial of the unknot is $-q^{\frac{1}{2}}-q^{-\frac{1}{2}}$.
The correct sign and power of $q^{\frac{1}{2}}$
is also explicitly specified in \cite{sBkyoto}.

This result is due to Lawrence,
who also achieved a similar result for the two-variable Jones polynomial
by a much more complicated construction \cite{rL98}.

%--------------------------------------------------
\section{The Hecke algebra}
\vskip-5mm \hspace{5mm }

We conclude with some speculation
about possible applications of the Blanchfield pairing
to the representation theory of Hecke algebras.
We first give a very brief overview of the basic theory of Hecke algebras.

Let $q \in \Complex \setminus \{0\}$.
The Hecke algebra $H_n(q)$, or simply $H_n$,
is the $\Complex$-algebra given by generators
$g_1,\dots,g_{n-1}$
and relations
\begin{itemize}
\item $g_i g_j = g_j g_i$ if $|i-j|>1$,
\item $g_i g_j g_i = g_i g_j g_i$ if $|i-j|=1$,
\item $(g_i - 1)(g_i + q) = 0$.
\end{itemize}
It is an $n!$-dimensional $\Complex$-algebra.
We are restricting to the ring $\Complex$ for convenience,
although other rings can be used.

Note that $H_n(1)$ is
the group algebra $\Complex S_n$ of the symmetric group $S_n$.
The Hecke algebra is called a ``quantum deformation'' of $\Complex S_n$.
The representation theory of $\Complex S_n$ is well understood
except when working over a field of finite characteristic in which $n! = 0$.
This is because the classical theory
sometimes requires one to divide by $n!$, the order of the group.
When studying $H_n$ it turns out to be useful to be able to divide by
$$(1+q+\dots+q^{n-1})(1+q+\dots+q^{n-2})\dots(1+q).$$
This is sometimes written as $[n]!$,
and can be thought of as a ``quantum deformation'' of $n!$.
Note that $[n]! = n!$ if $q=1$.
A {\em generic} value of $q$ is one for which $[n]! \neq 0$.
The non-generic values are
the primitive $k$th roots of unity for $k=2,\dots,n$.
The representation theory of $H_n$
is well understood for generic values of $q$,
but the non-generic values are the subject of ongoing active research.

One of the most important papers on this subject
is Dipper and James \cite{DJ86}.
For every partition $\lambda$ of $n$,
they define a $H_n$-module $S^\lambda$ called the {\em Specht module}.
They then define a bilinear pairing on $S^\lambda$,
which we denote $\langle \cdot,\cdot \rangle_{\rm DJ}$.
Let $S^\lambda_\perp$ denote the set of $u \in S^\lambda$
such that $\langle u,v \rangle_{\rm DJ} = 0$
for all $v \in S^\lambda$.
Let $D^\lambda$ be the quotient module $S^\lambda/S^\lambda_\perp$.
Dipper and James show that the non-zero $D^\lambda$
form a complete list of all distinct irreducible representations of $H_n$.
For generic values of $q$ we have $D^\lambda = S^\lambda$.
For non-generic values of $q$,
the $D^\lambda$ are not well understood.

Lawrence \cite{rL98} gave a homological definition of the Specht modules.
The construction begins with the action of $B_n$ on
a homology module of a configuration space.
The variable $t$ is then specialised to $-q^{-1}$,
and a certain quotient module is taken.
A detailed treatment of the case $\lambda = (n-2,2)$
is given in \cite{sBgeorgia}.

There is a Blanchfield pairing
on the Specht modules as defined by Lawrence.
It would be nice if this were the same
as the pairing defined by Dipper and James.
Unfortunately the Blanchfield pairing is sesquilinear,
whereas the pairing defined by Dipper and James is bilinear.
This problem can be overcome as follows.
Let $\rho \co D_n \to D_n$ be the conjugation map.
Let $\tilde{\rho}$ be the induced map on $H_k(\tilde{C},\tilde{B})$.
Then the pairing
$$\langle u,v \rangle' = \langle u,\rho(v) \rangle$$
can be shown to give a bilinear pairing on the Specht module.

For generic values of $q$,
this topologically defined pairing
is the same as the algebraically defined pairing of Dipper and James,
up to some renormalisation.
There is some evidence that
this can be made to work at non-generic values of $q$.
If so,
it would give rise to a new homological definition of the modules $D^\lambda$,
and new topological tools for studying them.
In any case,
it would be interesting to better understand
the behaviour of this Blanchfield pairing at roots of unity.

%--------------------------------------------------

\label{lastpage}

\end{document}